\documentclass{amsart}
\usepackage{amsopn}
\usepackage{amssymb, amscd}

\newcommand{\nc}{\newcommand}

\nc{\vg}{\mathfrak{v} }
\nc{\wg}{\mathfrak{w} }
\nc{\zg}{\mathfrak{z} }
\nc{\ngo}{\mathfrak{n} }
\nc{\kg}{\mathfrak{k} }
\nc{\mg}{\mathfrak{m} }
\nc{\bg}{\mathfrak{b} }
\nc{\ggo}{\mathfrak{g} }
\nc{\ggob}{\overline{\mathfrak{g}} }
\nc{\sog}{\mathfrak{so} }
\nc{\sug}{\mathfrak{su} }
\nc{\spg}{\mathfrak{sp} }
\nc{\slg}{\mathfrak{sl} }
\nc{\glg}{\mathfrak{gl} }
\nc{\cg}{\mathfrak{c} }
\nc{\rg}{\mathfrak{r} }
\nc{\hg}{\mathfrak{h} }
\nc{\tg}{\mathfrak{t} }
\nc{\ug}{\mathfrak{u} }
\nc{\dg}{\mathfrak{d} }
\nc{\ag}{\mathfrak{a} }
\nc{\pg}{\mathfrak{p} }
\nc{\sg}{\mathfrak{s} }
\nc{\pca}{\mathcal{P}}
\nc{\nca}{\mathcal{N}}
\nc{\lca}{\mathcal{L}}
\nc{\oca}{\mathcal{O}}
\nc{\mca}{\mathcal{M}}
\nc{\tca}{\mathcal{T}}
\nc{\aca}{\mathcal{A}}
\nc{\vp}{\varphi}
\nc{\ddt}{\frac{{\rm d}}{{\rm d}t}}
\nc{\im}{\mathtt{i}}

\nc{\SO}{{\mathrm SO}}
\nc{\Spe}{{\mathrm Sp}}
\nc{\Sl}{{\mathrm SL}}
\nc{\SU}{{\mathrm SU}}
\nc{\Or}{{\mathrm O}}
\nc{\U}{{\mathrm U}}
\nc{\Gl}{{\mathrm GL}}
\nc{\Se}{{\mathrm S}}
\nc{\Cl}{{\mathrm Cl}}
\nc{\Spin}{{\mathrm Spin}}
\nc{\Pin}{{\mathrm Pin}}

\nc{\RR}{{\Bbb R}}
\nc{\HH}{{\Bbb H}}
\nc{\CC}{{\Bbb C}}
\nc{\ZZ}{{\Bbb Z}}
\nc{\FF}{{\Bbb F}}
\nc{\NN}{{\Bbb N}}
\nc{\QQ}{{\Bbb Q}}
\nc{\PP}{{\Bbb P}}

\nc{\vs}{\vspace{.2cm}}
\nc{\ip}{\langle\cdot,\cdot\rangle}
\nc{\la}{\langle}
\nc{\ra}{\rangle}
\nc{\unm}{\frac{1}{2}}
\nc{\unc}{\frac{1}{4}}
\nc{\und}{\frac{1}{16}}
\nc{\no}{\vs\noindent}
\nc{\lam}{\Lambda^2\ggo^*\otimes\ggo}
\nc{\tang}{{\rm T}}
\nc{\dif}{{\rm d}}

\nc{\He}{\operatorname{Hess}}
\nc{\ad}{\operatorname{ad}}
\nc{\Ad}{\operatorname{Ad}}
\nc{\rank}{\operatorname{rank}}
\nc{\Irr}{\operatorname{Irr}}
\nc{\End}{\operatorname{End}}
\nc{\Aut}{\operatorname{Aut}}
\nc{\Inn}{\operatorname{Inn}}
\nc{\Der}{\operatorname{Der}}
\nc{\Ker}{\operatorname{Ker}}
\nc{\Iso}{\operatorname{I}}
\nc{\Diff}{\operatorname{D}}
\nc{\Lie}{\operatorname{L}}
\nc{\tr}{\operatorname{tr}}
\nc{\sen}{\operatorname{sen}}
\nc{\modu}{\operatorname{mod}}
\nc{\Ric}{\operatorname{Ric}}
\nc{\sym}{\operatorname{sym}}
\nc{\sca}{\operatorname{sc}}
\nc{\scalar}{{\sf s}}
\nc{\grad}{\operatorname{grad}}
\nc{\ricci}{\operatorname{ric}}
\nc{\Rin}{\operatorname{M}}
\nc{\Le}{\operatorname{L}}
\nc{\level}{\operatorname{level}}
\nc{\rad}{\operatorname{r}}
\nc{\abel}{\operatorname{ab}}

\newtheorem{theorem}{Theorem}[section]

\newtheorem{definition}[theorem]{Definition}
\newtheorem{remark}[theorem]{Remark}

\newtheorem{example}[theorem]{Example}

\title{Examples of Anosov diffeomorphisms}

\author{Jorge Lauret}

\address{Department of Mathematics, Yale University,
10 Hillhouse Box 208283 New Haven, CT 06520 USA (current affiliation). }
\email{jorge.lauret@yale.edu}

\address{FaMAF and CIEM, Universidad Nacional de C\'ordoba, 5000 C\'ordoba, Argentina}
\email{lauret@mate.uncor.edu}

\thanks{
Supported by CONICET and Guggenheim fellowships, and grants
from FONCyT and SeCyT UNC (Argentina).}

\begin{document}

\maketitle

\begin{abstract}
We prove that if $\ngo$ is any graded rational Lie algebra, then the simply connected
nilpotent Lie group $N\times N$ with Lie algebra $(\ngo\otimes\RR)\oplus(\ngo\otimes\RR)$
has a lattice $\Gamma$ such that the corresponding nilmanifold $(N\times N)/\Gamma$ admits an
Anosov diffeomorphism.  This gives a procedure to construct easily several explicit
examples of nilmanifolds admitting an Anosov diffeomorphism, and shows that a reasonable
classification up to homeomorphism (or even up to commensurability) of such
nilmanifolds would not be possible.
\end{abstract}

\section{Introduction}\label{intro}

Anosov diffeomorphisms play an important and beautiful role in dynamics as the notion
represents the most perfect kind of global hyperbolic behavior, giving examples of
structurally stable dynamical systems (see $\S$\ref{anonil} for a precise definition).
In \cite{Sml}, S. Smale raised the problem
of classifying the compact manifolds (up to homeomorphism) which admit an Anosov
diffeomorphism.  At this moment, the only known examples are of algebraic nature, namely
Anosov automorphisms of nilmanifolds and infranilmanifolds.  It is conjectured that
any Anosov diffeomorphism is
topologically conjugate to an Anosov automorphism of an infranilmanifold.  J. Franks
\cite{Frn} and A. Manning \cite{Mnn} proved the conjecture for Anosov diffeomorphisms
on infranilmanifolds themselves.

All this certainly highlights the problem of classifying all nilmanifolds which admit
Anosov automorphisms, which are easily seen in correspondence with hyperbolic automorphisms
of nilpotent Lie algebras over $\ZZ$.  Nevertheless, not too much is known on
the question since it is not so easy for an automorphism of a (real) nilpotent Lie
algebra being hyperbolic and unimodular at the same time.  We propose the following

\begin{definition}\label{anolie}
{\rm A rational Lie algebra $\ngo$ (i.e. with structure constants in $\QQ$) of
dimension $d$ is said to be {\it Anosov} if it admits
a hyperbolic automorphism $A$ (i.e. all their eigenvalues have absolute value
different from $1$)
such that $[A]_{\beta}\in\Gl(d,\ZZ)$ for some basis $\beta$ of $\ngo$, where $[A]_{\beta}$
denotes the matrix of $A$ with respect to $\beta$.  }
\end{definition}

It is well known that any Anosov Lie algebra is necessarily nilpotent, and it is easy
to see that the classification of nilmanifolds which admit an Anosov automorphism
is essentially equivalent to that of Anosov Lie algebras
(see $\S$\ref{anonil}).  If $\ngo$ is a rational Lie algebra, we call the real
Lie algebra $\ngo\otimes\RR$ the {\it real completion} of $\ngo$.

Curiously enough, the only explicit examples of Anosov Lie
algebras in the literature so far
are an infinite family of $6$-dimensional rational Lie algebras all
having real
completion isomorphic to $\hg_3\oplus\hg_3$ (see \cite{Mlf, Ito}) --where $\hg_3$ is the
3-dimensional
Heisenberg Lie algebra--, the free $k$-step nilpotent Lie algebras on $n$ generators
with $k<n$ (see \cite{Dn}, and also \cite{DkmMlf, Dkm} for a different approach) and
certain $k$-step nilpotent Lie algebras of dimension $d+\binom{d}{2}+...+\binom{d}{k}$
with $d\geq k^2$ (see \cite{Frd}).  For the known examples
of infranilmanifolds which are not nilmanifolds and admit Anosov automorphisms we refer
to \cite{Shb,Prt,Mlf2}.

This motivates the following natural questions:

\begin{itemize}
\item[(i)]  The minimal possible dimension for an Anosov non-abelian Lie
algebra is $6$ (see \cite{Mlf,Ito}).  Is there an Anosov Lie algebra (without a
nonzero abelian factor) in every dimension $d\geq 6$?.

\item[(ii)]  It is proved in \cite[Theorem 3.2]{Mlf} that any Anosov $k$-step nilpotent Lie
algebra $\ngo$ has $\dim{\ngo}\geq 2k+2$.  Is there for any $k\geq 2$ an Anosov $\ngo$ with
$\dim{\ngo}=2k+2$?.

\item[(iii)]  Does there exist a dimension $d$ for which there are infinitely many Anosov
Lie algebras such that their real completions are pairwise non-isomorphic?.

\item[(iv)]  Is it feasible a more or less explicit classification of all Anosov
Lie algebras (or at lest of their real completions) up to isomorphism?.

\item[(v)]  Does every Anosov automorphism on a nilmanifold come from the Anosov action
of some lattice in a semisimple Lie group? (see \cite{Dn,KS}).
\end{itemize}

The purpose of this paper is to give a procedure to construct explicit examples of Anosov
Lie algebras, which is inspired on the famous $6$-dimensional example given
by S. Smale in \cite{Sml} (due to A. Borel), but from the point of view of L. Auslander
and J. Scheuneman in \cite[Section 4]{AS}.  We prove that $\ngo\oplus\ngo$ is Anosov for any graded
rational nilpotent Lie algebra $\ngo$.  The construction is very direct and simple, but
it produces examples enough to answer positively questions (i) (for $d$ even), (ii)
and (iii), and negatively questions (iv) and (v).  Also, two examples of Anosov
Lie algebras which are not of the form $\ngo\oplus\ngo$ are given.  One of them has
dimension $9$, and so by considering direct sums we get examples for question (i) for
any $d\ne 7,11,13$.

We note that any Anosov Lie algebra is necessarily graded as they always admit semisimple
hyperbolic automorphisms (see \cite{AS}).

\vs
\noindent {\it Acknowledgements.}  I wish to express my gratitude to F. Grunewald for
pointing me out Remark \ref{gru}.  I am
also grateful to D. Fisher and R. Miatello for helpful observations.

\section{Anosov diffeomorphisms on nilmanifolds}\label{anonil}

A diffeomorphism $f$ of a compact differentiable manifold $M$ is called {\it Anosov}
if the tangent bundle $\tang M$ admits a continuous invariant splitting
$\tang M=E^+\oplus E^-$ such that $\dif f$ expands $E^+$ and contracts $E^-$ exponentially,
that is, there exist constants $0<c$ and $0<\lambda<1$ such that
$$
||\dif f^n(X)||\leq c\lambda^n||X||, \quad \forall X\in E^-, \qquad
||\dif f^n(Y)||\geq c\lambda^{-n}||Y||, \quad \forall Y\in E^+,
$$
for all $n\in\NN$.  The condition is independent of the Riemannian metric.

Let $N$ be a real simply connected nilpotent Lie group with Lie algebra $\ngo$.
Let $\vp$ be a hyperbolic automorphism of $N$, that is, all the eigenvalues
of its derivative $A=(\dif\vp)_e:\ngo\mapsto\ngo$ have absolute value different from $1$.
If $\vp(\Gamma)=\Gamma$ for some lattice $\Gamma$ of $N$ then $\vp$  defines an Anosov
diffeomorphism
on the nilmanifold $N/\Gamma$, which shall be called an {\it Anosov automorphism}.  The
subspaces $E^+$ and $E^-$ are obtained by left translation of the eigenspaces of eigenvalues
of $A$ of absolute value greater than $1$ and less than $1$, respectively.

We now review the well-known relationship with Definition \ref{anolie}
(see \cite{Dn,Ito,Dkm} for more detailed expositions).  The lattice $\Gamma$
contains, as a subgroup of finite index, some {\it full lattice}
$\Gamma_1$,
that is,
$\log(\Gamma_1)$ is a $\ZZ$-Lie subalgebra of $\ngo$ (i.e. with integer structure
constants).  Thus $\vp^m(\Gamma_1)
=\Gamma_1$ for some $m$ and hence $A^m$ is an hyperbolic
automorphism of the $\ZZ$-Lie algebra $\log(\Gamma_1)$.  This implies that the
rational completion $\log(\Gamma_1)\otimes\QQ$
is an Anosov Lie algebra.

Conversely, let $A$ be an automorphism of a rational Lie algebra $\ngo$ satisfying the
conditions in Definition \ref{anolie}.  Then
there exists a $\ZZ$-Lie subalgebra $\ngo_1$ of $\ngo$ such that
$\Gamma=\exp(\ngo_1)$ is a lattice in $N$, where $N$ is the simply connected nilpotent
Lie group
with Lie algebra $\ngo\otimes\RR$, the real completion of $\ngo$.  Now for some $m$,
$A^m(\ngo_1)=\ngo_1$ and hence $A^m$ determines an Anosov
automorphism of the nilmanifold $N/\Gamma$.  The second condition on $A$ in
Definition \ref{anolie} is equivalent to the fact that
the characteristic polynomial of $A$ has integer coefficients and constant term equal to
$\pm 1$ (see \cite{Dkm}).

Two lattices $\Gamma_1,\Gamma_2$ of $N$ are called {\it commensurable} if
$\Gamma_1\cap\Gamma_2$ has finite index in both $\Gamma_1$ and $\Gamma_2$, or equivalently,
the Lie algebras of their rational Mal'cev completions are
isomorphic as rational Lie algebras.  It is proved in \cite{Ito} and \cite{Dkm} that if
$\Gamma_1$ and $\Gamma_2$ are commensurable
then $N/\Gamma_1$ admits an Anosov automorphism if and only if $N/\Gamma_2$ does.  This
fact is the main reason for defining the notion of Anosov in the context of rational
Lie algebras.

Recall that if $\ngo_1$ and $\ngo_2$ are two non-isomorphic Anosov Lie algebras,
then the corresponding Anosov diffeomorphisms are not topologically conjugate since the
nilmanifolds $N_1/\Gamma_1$ and $N_2/\Gamma_2$ can never be homeomorphic.
Indeed, $\Pi_1(N_i/\Gamma_i)=\Gamma_i$ and $\ngo_i$ is the Lie algebra of the rational
Mal'cev completion of $\Gamma_i$.

Let $\ngo$ be a real nilpotent Lie algebra of dimension $d$.  We note that $\ngo$ has
a rational form which is Anosov (or equivalently, $N$ is the simply connected cover
of a nilmanifold admitting an Anosov automorphism) if and only if there exists a
$\ZZ$-basis $\beta$ of $\ngo$ and a hyperbolic $A\in\Aut(\ngo)$ such that
$[A]_{\beta}\in\Gl(d,\ZZ)$, where
$[A]_{\beta}$ is the matrix of $A$ in terms of $\beta$.

\section{Construction of examples}

A rational Lie algebra $\ngo$ is said to be {\it graded} if there exist (rational)
subspaces $\ngo_i$ of $\ngo$ such that
$$
\ngo=\ngo_1\oplus\ngo_2\oplus...\oplus\ngo_k \qquad \mbox{and} \qquad [\ngo_i,\ngo_j]\subset\ngo_{i+j}.
$$
Equivalently, $\ngo$ is graded when there are nonzero (rational) subspaces
$\ngo_{d_1},...,\ngo_{d_r}$, $d_1<...<d_r$,
such that $\ngo=\ngo_{d_1}\oplus...\oplus\ngo_{d_r}$ and if $0\ne [\ngo_{d_i},\ngo_{d_j}]$
then $d_i+d_j=d_k$ for some $k$ and $[\ngo_{d_i},\ngo_{d_j}]\subset n_{d_k}$.  Recall that
any graded Lie algebra is necessarily nilpotent.

\begin{theorem}\label{nn}
Let $\ngo$ be a graded rational nilpotent Lie algebra.  Then the direct sum
$\tilde{\ngo}=\ngo\oplus\ngo$ is Anosov.
\end{theorem}

\begin{proof}
Let $\{ X_1,...,X_d\}$ be a $\ZZ$-basis of $\ngo$ compatible with the gradation
$\ngo=\ngo_{d_1}\oplus...\oplus\ngo_{d_r}$, that is, each $X_i\in\ngo_{d_j}$ for
some $j$.  We denote by $\{ Y_i\}$ a copy of the basis $\{ X_i\}$, so that $\{ X_1,...,X_d,Y_1,...,Y_d\}$
is a basis of $\tilde{\ngo}$ and
\begin{equation}\label{basis}
[X_i,X_j]=\sum_{k=1}^{d} m_{ij}^kX_k, \qquad
[Y_i,Y_j]=\sum_{k=1}^{d} m_{ij}^kY_k, \qquad m_{ij}^k\in\ZZ.
\end{equation}
Every nonzero $\lambda\in\RR$ defines an automorphism $A_{\lambda}$ of $\ngo$ by
$$
A_{\lambda}|_{\ngo_{d_i}}=\lambda^{d_i}I,
$$
and also an automorphism $\tilde{A}_{\lambda}$ of $\tilde{\ngo}$ by
\begin{equation}\label{auto}
\tilde{A}_{\lambda}=
\left[\begin{array}{cc}
A_{\lambda}&\\
&A_{\lambda}^{-1}
\end{array}\right].
\end{equation}
If $a\in\ZZ$, $a\geq 2$, then the roots $\lambda,\lambda^{-1}$ of $x^2-2ax+1\in\ZZ[x]$
are
\begin{equation}\label{lambda}
\lambda=a+(a^2-1)^{\unm}, \qquad  \lambda^{-1}=a-(a^2-1)^{\unm},
\end{equation}
and hence $0<\lambda^{-1}<1<\lambda$.  Consider the new basis of $\tilde{\ngo}$ defined by
$$
\beta=\{ X_1+Y_1, (a^2-1)^{\unm}(X_1-Y_1),...,X_d+Y_d,(a^2-1)^{\unm}(X_d-Y_d)\}.
$$
It follows from (\ref{basis}) that
$$
\begin{array}{l}
[X_i+Y_i,X_j+Y_j]=\displaystyle{\sum_{k=1}^d} m_{ij}^k(X_k+Y_k),\\

[X_i+Y_i,(a^2-1)^{\unm}(X_j-Y_j)]=\displaystyle{\sum_{k=1}^d} m_{ij}^k(a^2-1)^{\unm}(X_k-Y_k), \\

[(a^2-1)^{\unm}(X_i-Y_i),(a^2-1)^{\unm}(X_j-Y_j)]=\displaystyle{\sum_{k=1}^d} (a^2-1)m_{ij}^k(X_k+Y_k).
\end{array}
$$
This implies that $\beta$ is also a $\ZZ$-basis.  On the other hand, it is easy to see
by using (\ref{lambda}) that if $T$ is the transformation whose matrix in terms of the
basis $\{ X_i,Y_i\}$ equals
$$
\left[\begin{array}{cc}
\lambda&0\\
0&\lambda^{-1}
\end{array}\right],
$$
then the matrix of $T$ with respect to the basis $\{ X_i+Y_i, (a^2-1)^{\unm}(X_i-Y_i)\}$
is given by
$$
B=\left[\begin{array}{cc}
a&a^2-1\\
1&a
\end{array}\right]\in\Sl(2,\ZZ).
$$
Therefore the matrix of $\tilde{A}_{\lambda}\in\Aut(\tilde{\ngo})$ (see
(\ref{auto})) in terms of the basis $\beta$ is given by
$$
\left[\tilde{A}_{\lambda}\right]_{\beta}=
\left[\begin{array}{ccccccc}
B^{d_1} &&&&&& \\
&\ddots&&&&&\\
&& B^{d_1}&&&&\\
&&&\ddots&&&\\
&&&&B^{d_r}&&\\
&&&&&\ddots&\\
&&&&&&B^{d_r}
\end{array}\right]\in\Sl(2d,\ZZ).
$$
Thus $\tilde{A}_{\lambda}$ determines a hyperbolic automorphism of $\tilde{\ngo}$
whose matrix in terms of $\beta$ is in $\Sl(2d,\ZZ)$, and hence $\tilde{\ngo}$ is
Anosov.
\end{proof}

\begin{remark}\label{irred}
{\rm In spite of the simply connected cover of the examples given by Theorem \ref{nn}
are all of the form $N\times N$, the nilmanifolds $(N\times N)/\Gamma$ which admit Anosov
automorphisms are not a direct product of nilmanifolds when $N$ is irreducible.  }
\end{remark}

We will now answer questions (i)-(v) in the
Introduction.  Theorem \ref{nn} will be applied in all cases without any further reference.

\no
{\bf Question (iv)}.  We first note that any two-step nilpotent Lie algebra is graded.
Therefore, the classification of all Anosov nilpotent Lie algebras contains the
classification of all rational two-step nilpotent Lie algebras.  This is equivalent to the
classification of all alternating bilinear maps
$$
\mu:\QQ^n\times\QQ^n\mapsto\QQ^m
$$
under the equivalence relation $\mu\simeq\lambda$ if and only if $\mu(gX,gY)=h\lambda(X,Y)$
for some $g\in\Gl(n,\QQ)$, $h\in\Gl(m,\QQ)$, which is considered for $m>2$ a wild
problem (see
\cite{GS}).  Moreover, the classification up to isomorphism of the real two-step nilpotent
Lie algebras admitting a rational form is also completely open.

\no
{\bf Question (iii)}.  We can answer this question positively for $d=14$ in the following
way.  For each $0<t<1$ consider the $7$-dimensional nilpotent Lie algebra $\ngo_t$ with
Lie bracket $\mu_t$ defined by
$$
\begin{array}{lcl}
\mu_t(X_1,X_2)=(1-t)^{\unm}X_3, && \mu_t(X_2,X_3)=X_5, \\

\mu_t(X_1,X_3)=X_4, && \mu_t(X_2,X_4)=X_6, \\

\mu_t(X_1,X_4)=t^{\unm}X_5, && \mu_t(X_2,X_5)=t^{\unm}X_7, \\

\mu_t(X_1,X_5)=X_6, && \mu_t(X_3,X_4)=(1-t)^{\unm}X_7. \\

\mu_t(X_1,X_6)=X_7, &&  \\
\end{array}
$$
Each Lie algebra $\mu_t$ is isomorphic to the Lie algebra
$\ggo_t$, denoted by $1,2,3,4,5,7_I:t$ in \cite[pp.494]{Sly}
(see also the curve $\tilde{g}(0,t,1,0,1,0,0,0)$ in \cite[5.2.3]{Mgn}).
The isomorphism $\vp_t:\ggo_t\mapsto\ngo_t$ is given by
$$
\begin{array}{l}
\vp_tX_1=X_1, \quad \vp_tX_2=t^{\unm}X_2, \quad \vp_tX_i=t^{\unm}(1-t)^{\unm}X_i,
\; i=3,4, \\ \\
\vp_tX_j=t(1-t)^{\unm}X_j, \; j=5,6,7.
\end{array}
$$
This proves that $\ngo_t$ ($0<t<1$) is a curve of pairwise non-isomorphic real Lie
algebras.  It is easy to check that
$\{ X_1,...,X_7\}$ is a $\QQ$-basis of $\ngo_{t_k}$ for every
$$
t_k=\frac{4k^2}{(k^2+1)^2}, \qquad k\in\NN,
$$
and that the subspaces $\ngo_i=\QQ X_i$ determines a gradation
$\ngo_{t_k}=\ngo_1\oplus...\oplus\ngo_7$ for any $k$.
Thus $\tilde{\ngo}_k=\ngo_{t_k}\oplus\ngo_{t_k}$ is Anosov for all $k\in\NN$,
and $\tilde{\ngo}_k\otimes\RR\equiv \tilde{\ngo}_{k'}\otimes\RR$ if and only if $k=k'$, since the function
$f(x)= \frac{4x^2}{(x^2+1)^2}$ is strictly decreasing in $[1,\infty)$.

We can also get examples for this question from any pairwise non-isomorphic infinite
sequence $\{ \ngo_k\}$ of real two-step nilpotent Lie algebras all admitting a rational
form.  There is an explicit example of such a sequence of dimension $10$ in \cite{Ggr},
and so $d=20$.

\no
{\bf Question (ii)}.  One can easily answer affirmatively this question by considering for
each $k\geq 2$ the following $k$-step nilpotent Lie algebra $\ngo$ of dimension $k+1$:
$$
[X_1,X_2]=X_3, \quad [X_1,X_3]=X_4, \quad ..., \quad [X_1,X_k]=X_{k+1},
$$
for which $\{ X_1,...,X_{k+1}\}$ is clearly a $\ZZ$-basis, and $\ngo$ is graded by
$\ngo=\QQ X_1\oplus...\oplus\QQ X_{k+1}$.  Thus $\tilde{\ngo}=\ngo\oplus\ngo$ is Anosov,
$(2k+2)$-dimensional and $k$-step nilpotent.  Actually, any $d$-dimensional graded
rational Lie algebra which is filiform (i.e. $(d-1)$-step nilpotent) provides an example
as required in this question.

\no
{\bf Question (v)}.  It is easy to prove that $\Aut(\ngo\oplus\ngo)$ is solvable for any
$\ngo$ defined in question (ii) with $k\geq 3$ (and for any filiform example with $d\geq 4$
as well).  Thus $\Aut(\ngo\oplus\ngo)$ does not contain any semisimple Lie subgroup,
which implies that the Anosov automorphism on $(N\times N)/\Gamma$ determined by
$\tilde{A}_{\lambda}$ (see (\ref{auto})) does not come from the Anosov action of any
lattice in a semisimple Lie group (compare with \cite{Dn}).

\no
{\bf Question (i)}.  We have required examples without a nonzero
abelian factor in order to avoid Lie algebras of the form
$\ngo\oplus\QQ^r$ with $\ngo$ Anosov.  For $r\geq 2$, the corresponding nilmanifolds
$N/\Gamma\times T^r$, where $T^r$ is
the $r$-dimensional torus, clearly admit Anosov automorphisms. The examples considered
in question (ii)
also give a positive answer to this question for any dimension $d\geq 6$ which is
even.  Thus, by using the $9$-dimensional Anosov Lie algebra given in Example \ref{libre},
we can consider direct sums and get examples of Anosov Lie algebras without an abelian
factor for any dimension $d\geq 15$.  Thus the question only remains open for
$d=7,11,13$.  It should be noted that the examples of odd dimensions $d\geq 15$ are all
direct products of nilmanifolds (compare with Remark \ref{irred}).

\vspace{.2cm}

We finally want to give two examples of Anosov two-step nilpotent Lie algebras which are
not of the form $\ngo\oplus\ngo$.  We note that the construction of the first one is
definitely similar in spirit to that in the proof of Theorem \ref{nn}, from where
it can be deduced that such
a construction might be just a particular case of a much more general procedure.

\begin{example}
{\rm Let $\ngo$ be the two-step nilpotent Lie algebra with
basis
$$
\{ X_1,X_2,X_3,Y_1,Y_2,Y_3,Z_1,Z_2\}
$$
and Lie bracket defined by
$$
[X_1,X_2]=Z_1, \quad [X_1,X_3]=Z_2, \quad [Y_1,Y_2]=Z_1, \quad [Y_1,Y_2]=Z_2.
$$
Consider $A\in\Aut(\ngo)$ given by
$$
\begin{array}{l}
AX_1=\lambda X_1, \quad   AX_2=\lambda X_2, \quad AX_3=\lambda^{-3} X_3, \\ \\

AY_1=\lambda^{-1} Y_1, \quad AY_2=\lambda^3 Y_2, \quad AY_3=\lambda^{-1} Y_3, \\ \\

AZ_1=\lambda^2 Z_1, \quad AZ_2=\lambda^{-2} Z_2,
\end{array}
$$
where $\lambda$ is as in (\ref{lambda}).  It is easy to see that
$$
\begin{array}{rl}
\beta=&\left\{X_1+Y_1, (a^2-1)^{\unm}(X_1-Y_1),X_2+Y_3,(a^2-1)^{\unm}(X_2-Y_3),X_3+Y_2,\right. \\
& \left.(a^2-1)^{\unm}(X_3-Y_2),Z_1+Z_2,(a^2-1)^{\unm}(Z_1-Z_2)\right\}
\end{array}
$$
is a $\ZZ$-basis of $\ngo$ and that if $B=
\left[\begin{array}{cc}
a&a^2-1\\
1&a
\end{array}\right]$ then
$$
\left[A\right]_{\beta}= \left[\begin{array}{cccc}
B&&& \\

&B&&\\

&&B^{-3}&\\
&&&B^2
\end{array}\right]\in\Sl(8,\ZZ),
$$
showing that $\ngo$ is Anosov.      }
\end{example}

The following example is just a modification of the free $2$-step nilpotent Lie algebra
on $3$ generators (see \cite{DkmMlf}).

\begin{example}\label{libre}
{\rm Consider the $(3r+3)$-dimensional $\ZZ$-Lie algebra
$$
\ngo=\ZZ^3\oplus...\oplus\ZZ^3\oplus\Lambda^2\ZZ^3
$$
with Lie bracket defined by
$$
[v_1+...+v_r,w_1+...+w_r]=v_1\wedge w_1+...+v_r\wedge w_r.
$$
If $A\in\Gl(3,\ZZ)$ has eigenvalues $\lambda_1,\lambda_2,\lambda_3$ and
their products $\lambda_i\lambda_j$ all of absolute value different from $1$ then it
easy to see that
$$
\left[\begin{array}{cccc}
A&&& \\

&\ddots&&\\

&&A&\\
&&&\Lambda^2A
\end{array}\right]\in\Gl(3r+3,\ZZ),
$$
is an hyperbolic automorphism of $\ngo$.  For instance, one can take
$$
A= \left[\begin{array}{ccc}
1&1&1 \\

1&2&2\\

1&2&3
\end{array}\right]
$$
(see \cite{DkmMlf}).  Thus $\ngo\otimes\QQ$ is Anosov of dimension $3r+3$. }
\end{example}

\begin{remark}\label{gru}
{\rm It was pointed me out by F. Grunewald that the construction given in Theorem \ref{nn}
can be generalized by using the ring of integers in any totally real number field $F$
over $\QQ$ (we have used $F=\QQ\oplus (a^2-1)^{\unm}\QQ$).  The proof is essentially the
same, it only uses Dirichlet Unit Theorem (compare with \cite[Proposition 3.2]{Ito} and
\cite{Frd}).  If the degree of $F$ is $m$ then we obtain Anosov Lie algebras of the form
$\ngo\oplus...\oplus\ngo$ ($m$ times) for any graded rational nilpotent Lie algebra
$\ngo$. }
\end{remark}

\end{document}